\documentclass{amsart}

\usepackage{amssymb}
\usepackage{a4wide}
\usepackage{graphicx}

\newcommand{\arr}{\longrightarrow}
\newcommand{\aut}{{\operatorname{aut}}}
\newcommand{\autxs}{{\aut(\xs)}}
\newcommand{\eker}{{\mathcal E}}
\newcommand{\kn}[1]{\mathfrak{K}(#1)}
\newcommand{\kv}{\kn{v}}
\newcommand{\kwv}{\kn{w,v}}
\newcommand{\img}[1]{\mathop{\mathrm{IMG}}(#1)}
\newcommand{\M}{{\mathcal M}}
\newcommand{\R}{{\mathbb R}}
\newcommand{\rel}{{\mathcal R}}
\newcommand{\symm}{{\mathfrak{S}_2}}
\newcommand{\T}{\R/\Z}
\newcommand{\xs}{X^*}
\newcommand{\Z}{{\mathbb Z}}
\newcommand{\N}{{\mathbb N}}
\newcommand{\ZZ}[1]{\frac{\Z}{#1\Z}}
\newcommand{\nuke}{\mathcal{N}}
\newcommand{\C}{\mathbb{C}}

\newbox\pairbox
\def\pair<#1,#2>{{\mathsurround=0pt
    \setbox\pairbox\hbox{$\left\langle#1,\;#2\right\rangle$}
    \left\langle\kern-0.35\ht\pairbox
      \copy\pairbox\kern-0.35\ht\pairbox\right\rangle}}

\newtheorem{theorem}{Theorem}[section]
\newtheorem{corollary}[theorem]{Corollary}
\newtheorem{proposition}[theorem]{Proposition}
\newtheorem{lemma}[theorem]{Lemma}

\theoremstyle{definition}
\newtheorem{defi}{Definition}

\begin{document}
\title{Iterated Monodromy Groups of Quadratic Polynomials, I}
\author{Laurent Bartholdi}
\author{Volodymyr V.~Nekrashevych}
\date{7 November 2006}
\begin{abstract}
  We describe the iterated monodromy groups associated with
  post-critically finite quadratic polynomials, and make explicit their
  connection to the `kneading sequence' of the polynomial.

  We then give recursive presentations by generators and relations for
  these groups, and study some of their properties, like torsion and
  `branchness'.
\end{abstract}
\maketitle

\section{Introduction}
Symbolic dynamics of quadratic complex polynomials are
traditionally studied by their ``kneading sequence'', an infinite
sequence over the alphabet $\{0,1,*\}$ that encodes symbolically
the dynamics of the map (see~\cite{henkdierk,keller}). This
sequence is (pre)periodic if the polynomial has finite
post-critical set, and we shall make that assumption here.

A construction by the second author associates a finitely generated
group, defined by automata, with such a polynomial. This group is the
iterated monodromy group of the polynomial, and acts on a binary
rooted tree.

In this paper, we show that the automata defining the group may be
chosen in a particularly simple manner. These automata are
``bounded'', i.e.\ their activity is trivial away from a ray in
the tree. If the tree's rays are labeled by infinite $\{0,1\}$
sequences, then the rays on which activity is non-trivial are the
translates of the periodic part of the kneading sequence of the
polynomial.

Actually, we study a class of groups containing all the iterated
monodromy groups of post-critically finite quadratic polynomials as a
proper subset. The first set of such groups $\kv$ corresponds to
kneading sequences of the form $(v*)^\omega$, and the second set of
groups $\kwv$ corresponds to kneading sequences of the form
$w(v^\omega)$. Since not every sequence of this form is realizable as
a kneading sequence, not all groups $\kv$ and $\kwv$ are iterated
monodromy groups of quadratic polynomials. In particular, if $v$ is
periodic (is a proper power of a word), then $\kwv$ is not an iterated
monodromy group of any polynomial (they are iterated monodromy groups
of \emph{obstructed topological polynomials}). In all the other cases
the corresponding groups are the iterated monodromy groups of some
post-critically finite polynomials of degree $2^n$. For more details
see a general description of the iterated monodromy groups of
post-critically finite polynomials in~\cite{nek:book}.

We describe elementary properties of the groups $\kv$ and $\kwv$.
Unsurprisingly, the groups are of a very different nature
depending on whether the kneading sequence is periodic (for $\kv$)
or pre-periodic (for $\kwv$):
\begin{itemize}
\item if the sequence is periodic, these groups are weakly branch, but
  not branch; they are torsion-free;
\item if the sequence is pre-periodic, these groups are branch, and
  contain elements of arbitrarily high 2-power order, as well as
  elements of infinite order.
\end{itemize}
We give in all cases recursive presentations for these groups.

The paper is organised as follows: \S\ref{ss:img} contains necessary
recollections on iterated monodromy groups; \S\ref{ss:kv} describes
the groups $\kv$ associated with periodic kneading sequences;
\S\ref{ss:kwv} describes the groups $\kwv$ associated with pre-periodic
kneading sequences; and \S\ref{ss:kneading} proves that these groups
are indeed the iterated monodromy groups of polynomials with
prescribed kneading sequence.

Our aim in this paper was to show that kneading sequences and
bounded automata are really equivalent descriptions for
post-critically finite quadratic polynomials and to describe the
basic algebraic properties of the iterated monodromy groups.

Some of the groups we consider already appeared in the literature, and
the present paper extends or unifies algebraic results concerning
them. The groups $\kv$ and $\kwv$ that coincide with previously known
ones are
\begin{itemize}
\item $\kn{}=\Z$;
\item $\kn{0,1}=\kn{1,0}=$ the infinite dihedral group;
\item $\kn{0}$, an ``amenable non-subexponentially-amenable'' group
  considered in~\cite{barthvirag}, where it is called the ``Basilica
  group'';
\item $\kn{0,11}$, a group whose growth is studied
  in~\cite{ershler:growth};
\item $\kn{0,111}$, a group whose Lie algebra is studied
  in~\cite{bartholdi-g:lie}, and whose presentation and growth are
  studied in~\cite{bartholdi-g:parabolic}, where it is called the
  ``Grigorchuk overgroup''. This group contains Grigorchuk's example
  of an infinite torsion group~\cite{grigorchuk:80_en} as a
  subgroup.
\end{itemize}

\subsection{Notation} The notation we use is quite standard, with
groups acting on the right. We therefore write $[a,b]=a^{-1}b^{-1}ab$,
and $a^b=b^{-1}ab$. More generally, we write
$g^{k_1h_1+\dots+k_nh_n}=(g^{k_1})^{h_1}\dots (g^{k_n})^{h_n}$ for
$k_i\in\Z$ and $g,h_i$ group elements. The derived subgroup $[G,G]$ of
$G$ is written $G'$. The cyclic group of order $n$ is written $\ZZ n$.

The notation $a=\pair<a_0,a_1>\sigma$ describes an automorphism $a$ of
the binary tree. It first acts as $a_0$ on the left subtree and as
$a_1$ on the right subtree, and then applies the permutation $\sigma$
to these two subtrees.

On the other hand, if $g$ is an element of a free group $F$ that acts
(in principle non-faithfully) on a tree, this action will in our cases
be described by a map $\Psi:F\to(F\times F)\rtimes2$. Then the
notation $\Psi(g)=\pair<g_0,g_1>\sigma$ describes the image of $g$ in
using the same convention as above.

If $v$ is a finite word over an alphabet $X$ then $|v|$ denotes
its length, i.e., such $n$ that $v\in X^n$.

\section{Iterated monodromy and self-similar groups}\label{ss:img}

\subsection{Iterated monodromy groups}
We present in this subsection a review of the notions and results
related to iterated monodromy groups. More details can be found
in~\cite{nek:book}.

A polynomial $f\in\C[z]$ is said to be \emph{post-critically
finite} if the orbit of every critical point is finite. The
union $P_f=\bigcup_{n\ge 1}f^{\circ n}(C_f)$ of the orbits of
critical points is called the \emph{post-critical set} of $f$;
here $C_f$ is the set of critical points of $f$.

Suppose that $f$ is a post-critically finite polynomial. We
consider it as a branched self-covering of the complex plane.
Moreover, it is a covering map $f:\C\setminus
f^{-1}(P_f)\arr\C\setminus P_f$ of the punctured plane
$\C\setminus P_f$ by its open subset $\C\setminus f^{-1}(P_f)$.

Consider now the general situation of a covering
$f:\mathcal{M}_1\arr\M$ of a path connected and locally path
connected topological space $\M$ by an open subset $\M_1$. We can
iterate the partially define map $f$ to obtain coverings
$f^n:\M_n\arr\M$ of $\M$ by open subsets.

Take a basepoint $t\in\M$. Then the fundamental group $\pi_1(\M,
t)$ acts naturally on the set of preimages $f^{-n}(t)$: the image
of a point $z\in f^{-n}(t)$ under the action of a loop
$\gamma\in\pi_1(\M, t)$ is equal to the endpoint of the unique
$f^n$-preimage of $\gamma$ that starts at $z$. Such a preimage
exists and is unique, since $f^n:\M_n\arr\M$ is a covering.

In this way, we get an action of $\pi_1(\M, t)$ on the formal disjoint
union $\sqcup_{n\ge 0}f^{-n}(t)$, the \emph{backward orbit} of $t$.
The quotient of $\pi_1(\M, t)$ by the kernel of its action is called the
\emph{iterated monodromy group} of $f$ and is denoted $\img{f}$.

The backward orbit $T$ has a natural structure of a rooted tree.
The root is the unique element $t$ of the set $f^{-0}(t)$ and a
vertex $z\in f^{-n}(t)$ of the $n$th level of the tree is
connected to the vertex $f(z)\in f^{-(n-1)}(t)$ of the $(n-1)$-st
level.

It is easy to see that the actions of the fundamental group and of
the iterated monodromy group preserve the structure of the rooted
tree.

Suppose that the degree of the covering $f$ (i.e., the number of
preimages of a point) is finite and is equal to $d$. Choose an
\emph{alphabet} $X$ of $d$ letters and consider the \emph{tree of
words} $\xs$ over this alphabet, i.e., the free monoid generated
by $X$, in which a word $v$ is connected to all words of the form
$vx$ for $x\in X$. The root of $\xs$ is the empty word. Since both
the tree of words $\xs$ and the tree of preimages $T$ are regular
$d$-trees, they are isomorphic.

A particularly nice class of isomorphisms $\Lambda:\xs\arr T$ are
constructed in the following way. Choose an arbitrary bijection
$\Lambda:X\arr f^{-1}(t)$ of the first level of the tree $\xs$
with the first level of the tree $T$. Choose also a collection of
\emph{connecting paths} $\ell_x$ from $t$ to $\Lambda(x)$ in $\M$.

We define now the isomorphism $\Lambda$ inductively
level-by-level. It is already defined on the zeroth and first
levels. Suppose that it is defined on the $n$th level. Let $v\in
X^n$ be an arbitrary vertex of the $n$th level and let $x\in X$ be
an arbitrary letter. Then the path $\ell_x$ has a unique
$f^n$-preimage starting at $\Lambda(v)$. We declare $\Lambda(xv)$
to be the endpoint of this preimage.

The map $\Lambda:\xs\arr T$ that we just constructed is an isomorphism
of rooted trees, and the following proposition makes it possible to
compute the action of elements of the fundamental group on the tree
$T$.

\begin{proposition}
\label{pr:imgformula} Let us conjugate the action of $\pi_1(\M,
t)$ on the tree $T$ by the isomorphism $\Lambda:\xs\arr T$
constructed above. Then the resulting action of $\pi_1(\M, t)$ on
$\xs$ is computed by the following recursive formula:
\[\gamma(xv)=y(\ell_x\gamma_x\ell_y^{-1})(v),\]
where $x\in X$, $v\in\xs$ and $\gamma\in\pi_1(\M, t)$ are
arbitrary, $\gamma_x$ is the $f$-preimage of $\gamma$ starting at
$\Lambda(x)$ and $y$ is such that $\Lambda(y)$ is the end of
$\gamma_x$.
\end{proposition}

\subsection{Automata}
Recursive formul\ae\ as in Proposition~\ref{pr:imgformula} are
conveniently interpreted in terms of automata theory. We interpret
$\img{f}$ (or $\pi_1(\M, t)$) as an automaton which being in a state
$\gamma$ and reading as input a letter $x$ produces as output the
letter $y$ and moves to state $\ell_x\gamma_x\ell_y^{-1}$. Then it is
ready to read a new input.

We describe such automata (or their subsets) by their \emph{Moore
diagram}. It is a graph whose vertices are identified with the
states of the automaton in which we have an arrow from a state
$q_1$ to a state $q_2$ labeled by a letter $x\in X$ if the
automaton being in the state $q_1$ after reading the letter $x$
goes to the state $q_2$. We also label the vertices of the
automaton by the transformations of the letters that they define.

In particular, if $X=\{0, 1\}$ and the automaton describes an action
of a group $G$ on the tree $\xs$, then every state $g\in G$ acts
either by the rule $g(0w)=0g_0(w)$ and $g(1w)=1g_1(w)$ or by the rule
$g(0w)=1g_0(w)$ and $g(1w)=0g_1(w)$, where $g_0, g_1$ are next states
of the automaton and $w\in\xs$ is arbitrary. In the first case we say
that $g$ is \emph{inactive} and write $g=\pair<g_0, g_1>$, in the
second case we say that $g$ is \emph{active} and write $g=\pair<g_0,
g_1>\sigma$. In both cases we have in the Moore diagram an arrow from
$g$ to $g_i$ labeled by $i$. We will mark the inactive states in a
Moore diagram by white dots, and the active states by black dots.

The notations $g=\pair<g_0, g_1>$ and $g=\pair<g_0, g_1>\sigma$ come
from the decomposition of the automorphism group $\autxs$ of the
binary rooted tree into the (permutational) wreath product
$\autxs\wr\symm=\autxs^X\rtimes\symm$. The factor $\symm$ acts on the
tree just by its action on the first letter of words, and the factor
$\autxs^X$ acts separately on the each of the subtrees $x\xs$:
\[(0v)^{\pair<g_0, g_1>}=0(v^{g_0}),\quad (1v)^{\pair<g_0, g_1>}=1(v^{g_1}).\]
Then $\sigma\in\symm$ denotes the transposition and $\pair<g_0,
g_1>\sigma$ is equal to the composition of $\pair<g_0, g_1>$ and
$\sigma$.

\subsection{Self-similar groups}

\begin{defi}
A group $G$ acting faithfully on the tree $\xs$ is said to be
\emph{self-similar} if for every $g\in G$ and every $x\in X$ there exist
$h\in G$ and $y\in X$ such that \[g(xw)=yh(w)\] for all $w\in\xs$.
\end{defi}

In other terms a self-similar group is an automaton such that the
set of transformations of $\xs$ defined by its states is a group
with respect to composition.

In particular, Proposition~\ref{pr:imgformula} shows that iterated
monodromy groups act on $\xs$ self-similarly (i.e., they are
self-similar groups).

If $G$ is a self-similar group acting on the binary tree then
each of its elements $g$ is decomposed either as $g=\pair<g_0, g_1>$
or as $g=\pair<g_0, g_1>\sigma$ according to its action on $X\subset\xs$
and, by the definition of self-similarity, we have $g_0, g_1\in
G$.

Hence every self-similar group $G$ comes with an associated
\emph{wreath recursion} $G\arr G\wr\mathfrak{S}_X$. On the other hand,
every such a homomorphism defines recursively an action of $G$ on
$\xs$ (which is also called self-similar, though it may be
non-faithful).

If $\Psi:G\arr G\wr\mathfrak{S}_X$ is a homomorphism, then for $g\in G$ and
$v\in\xs$ we define the \emph{restriction} $g|_v$ recursively by
the rule
\[g=\pair<g|_0, g|_1>\text{\ or\ }g=\pair<g|_0, g|_1>\sigma\]
and
\[g|_{vx}=g|_v|_x\]
for all $v\in\xs$ and $x\in X$. We also assume that
$g|_\emptyset=g$, where $\emptyset$ denotes the empty word.

If the action of $G$ on $\xs$ is faithful, then the restrictions
are uniquely defined by the condition
\[g(vu)=g(v)g|_v(u)\]
for $v, u\in\xs$.

\begin{proposition}\label{prop:torsionfree}
  Let $G$ be a self-similar group and let $N$ be a normal subgroup of
  $G$ which belongs to the stabilizer of the first level. If $G/N$ is
  torsion free, then so is $G$.
\end{proposition}
\begin{proof}
  Suppose that $g\in G$ has finite order. Since $G/N$ is torsion free,
  we have $g\in N$. This implies that $g$ belongs to the stabilizer of
  the first level, so $(g^n)|_x=(g|_x)^n$ for all $x\in X$. The
  elements $g|_x$ therefore also have finite order. Repeating the same
  argument, we obtain inductively that $g$ fixes $\xs$, so $g=1$.
\end{proof}

\begin{defi}
A self-similar group $G$ is called \emph{recurrent} (or
\emph{self-replicating}) if its action is transitive on the first
level of the tree $\xs$ and if for some (and thus for all) $x\in
X$ the homomorphism $g\mapsto g|_x:G_x\arr G$ is onto, where $G_x$
is the stabilizer $G_x$ of $x$ in $G$.
\end{defi}

It is easy to prove that if a self-similar group is recurrent, then it
is transitive on every level of the tree $\xs$ (the group is then
called \emph{level-transitive}).

\begin{defi}
A homomorphism $\Psi:G\arr G\wr\mathfrak{S}_X$ is \emph{contracting} if
there exits a finite set $N\subset G$ such that for every $g\in G$
there exists $n\in\N$ such that $g|_v\in N$ for all words
$v$ of length $\ge n$. The smallest such set $N$ is called the
\emph{nucleus} of $\Psi$.
\end{defi}

We say that a self-similar group is contracting if the associated
wreath recursion is contracting. The nucleus of the recursion is
called then the nucleus of the group.

It follows directly from the definitions that the nucleus is an
automaton, i.e., that for every element $g$ of the nucleus and for
every $x\in X$ the restriction $g|_x$ is in the nucleus.

\subsection{Bounded automata}
\begin{defi}
An automorphism $g\in\autxs$ is \emph{finite state} if it is
defined by a finite automaton, or, equivalently, if the set
$\{g|_v\;:\;v\in\xs\}$ is finite.

An automorphism $g\in\autxs$ is called \emph{bounded} if it is
finite-state and the sequence
\[q_n=|\{v\in X^n\;:\;g|_v\ne 1\}|\]
is bounded. Here $1$ is the identity tree automorphism.
\end{defi}

It is not hard to see that the set of all bounded automorphisms of
the tree $\xs$ is a group.

Bounded automata and tree automorphisms were defined and studied
for the first time by S.~Sidki in~\cite{sid:cycl}. The following
description of bounded automata follows from his results.

\begin{defi}
An automorphism $g\in\autxs$ is called \emph{finitary} if there
exists $n\in\N$ such that $g|_v=1$ for all $v\in\xs$ of
length more than $n$.

An automorphism $g\in\autxs$ is called \emph{directed} if there
exists $v\in\xs$ such that $g|_v=g$ and $g|_u$ is finitary for
every $u\in\xs$ such that $u\ne v$ and $|u|=|v|$. Then the
infinite word $v^\omega=vvv\ldots$ is called the \emph{kneading
sequence} of $g$.
\end{defi}

If $g$ is finitary, then it acts only on the first $n$ letters of
every word for some fixed $n$. If it is directed, then all of its
non-trivial action on the tree is concentrated around the path
described by the kneading sequence.

It is easy to see that every finitary and every directed
automorphism is bounded.

\begin{proposition}
If $g\in\autxs$ is bounded, then there exists $n$ such that $g|_v$
for every $v\in X^n$ is either finitary or directed.
\end{proposition}

We have the following properties of groups of bounded automata.

\begin{theorem}[\cite{bondnek}]
If $G$ acts on $\xs$ by bounded automata, is finitely generated
and self-similar, then it is contracting.
\end{theorem}

\begin{theorem}[\cite{bknv}] The group of all bounded
automorphisms of the tree $\xs$ is amenable.
\end{theorem}

\subsection{Branch groups}

We say that a group $G$ acting on a regular tree $X^*$ is
\emph{regular weakly branch on $H\le G$} if $H$ is non-trivial and
contains the geometric direct product $H\times\dots\times H$ with
$|X|$ factors. Here by \emph{geometric direct product} we mean the
group generated by copies of $G$ acting disjointly on all subtrees
$xX^*$ with $x$ ranging over $X$. The group $G$ is \emph{regular
  branch on $H$} if it is weakly branch on some subgroup $H$ that has
finite index in $G$.

\def\0{\cite[Theorem 6.9]{handbook:branch}}
\begin{proposition}[\0]\label{prop:containsW}
  Let $G$ be a regular weakly branch group on its subgroup $L$, and
  suppose that a subgroup $K$ of $L$ has a regular orbit on $X^n$ for
  some $n$. Then $G$ contains $\lim \wr^nK=\bigcup_{n\ge0}K\wr\dots\wr
  K$.

  In particular, if there exists such a $K\cong\ZZ p$, then $G$
  contains every finite $p$-group.
\end{proposition}
\begin{proof}
  Since $G$ is weakly branch, it contains the subgroup $K_0=K$, and
  for all $i\ge1$ the subgroup $K_i=1\times\dots\times
  K_{i-1}\times\dots\times1$ with $|X|^n$ factors, where the
  $K_{i-1}$ is at any position on a regular orbit of $K$.

  The group generated by $K_0\cup K_1\cup\dots$ is isomorphic to
  $\lim\wr^n K$.

  It is known~\cite{kaloujnine-krasner} that every extension embeds in
  a wreath product; since finite $p$-groups are iterated extensions of
  $\ZZ p$, the second claim follows.
\end{proof}

\section{The groups $\kv$}\label{ss:kv}
Let $v=x_1x_2\ldots x_{n-1}$ be a word over the alphabet $X=\{0,1\}$.
We denote by $\kv$ the subgroup of $\autxs$ generated by the elements
$a_1, \ldots, a_n$ defined by
\[a_1=\pair<1, a_n>\sigma,\qquad
a_{i+1}=\begin{cases}\pair<a_i, 1> & \text{ if }x_i=0,\\
  \pair<1, a_i> & \text{ if }x_i=1,\end{cases}\text{ when }1\le i<n.
\]

\begin{figure}
\begin{center}
  \includegraphics{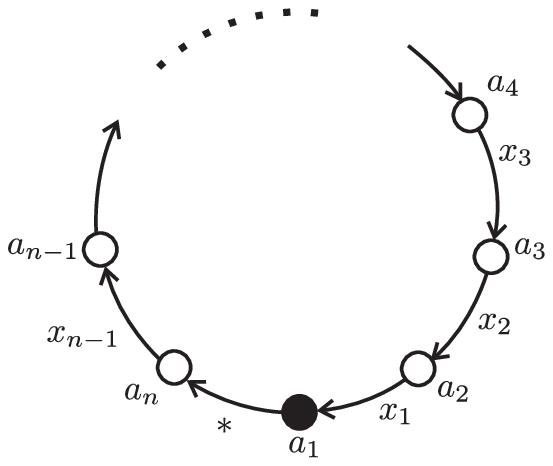}
  \caption{Automaton generating $\kv$}\label{fig:kv}
  \end{center}
\end{figure}

In other words, the group $\kv$ is generated by the automaton whose
Moore diagram is shown on Figure~\ref{fig:kv}.  In this diagram, only
the edges leading to non-trivial states are drawn, and the active
state $a_1$ is labelled by $\sigma$.

Note that, in the group $\kv$, if we replace the generators $a_i$
by their inverses $a_i'$, we obtain the description
\[a_1'=\pair<a_n',1>\sigma,\qquad
a_{i+1}'=\begin{cases}\pair<a_i', 1> & \text{ if }x_i=0,\\
  \pair<1, a_i'> & \text{ if }x_i=1,\end{cases}\text{ when }1\le i<n.
\]
Therefore, if we modify the definition of $\kv$ by setting
$a_1=\pair<a_n, 1>\sigma$, we do not change the group.

Let us denote by $v'$ the word obtained from the word $v$ by changing
$0$ to $1$ and $1$ to $0$. It is easy to see then that the groups
$\kv$ and $\kn{v'}$ are conjugate; more precisely,
$\kn{v'}=\kv^\alpha$, where $\alpha\in\autxs$ is defined by the
recursion $\alpha=\pair<\alpha,\alpha>\sigma$. Note that $\alpha$ also
interchanges $0$ and $1$, i.e.\ $v^\alpha=v'$.

\begin{lemma}
  The group $\kv$ is recurrent and level-transitive.
\end{lemma}
\begin{proof}
  The projections of the $a_{i+1}$, $a_{i+1}^{a_1}$ and $a_1^2$, all
  fixing the first level of $\xs$, give all generators; the action is
  therefore recurrent.

  To check that the action of a recurrent group is level-transitive,
  it suffices to check that it is transitive on the first level; this
  is achieved by $a_1$.
\end{proof}

\subsection{Wreath recursion}
Let us denote by $F_n$ the free group on $\{a_1,\dots,a_n\}$, and by
$\symm=\{1, \sigma\}$ the symmetric group on $X=\{0, 1\}$. Let
$\Psi:F_n\arr F_n\wr\symm$ be given by the recursive definition of the
group $\kv$, i.e.,
\begin{equation}\label{eq:kv:Psi}
\Psi(a_1)=\pair<1, a_n>\sigma,\qquad
\Psi(a_{i+1})=\begin{cases}\pair<a_i, 1> & \text{ if }x_i=0,\\
  \pair<1, a_i> & \text{ if }x_i=1,\end{cases}\text{ when }1\le i<n.
\end{equation}

We can iterate the map $\Psi$ to obtain a homomorphism
$\Psi^k:F_n\arr F_n\wr\symm\wr\cdots\wr\symm$, where $\symm$
appears $k$ times in the iterated permutational wreath product.

\begin{lemma}\label{lem:kv:1}
  The wreath recursion $\Psi$ is contracting:
  if $\Psi^n(g)=\pair<g_0,\dots,g_{2^n-1}>\pi$ and $\|g\|>2$, then
  $\|g_i\|<\|g\|$ for all $i$.

  The nucleus of $\kv$ (and of the wreath recursion $\Psi$)
  may be expressed as follows: let $d$ be maximal
  such that $vx=u^d$ for some $x\in\{0,1\}$. Then the nucleus of $\kv$
  contains $1+n(d+1)$ elements. Taking indices modulo $n$, they are
  $1,a_i^{\pm1}$, and $a_i^{\epsilon}a_{i+(n/d)j}^{-\epsilon}$ for any
  $j\in\{1,\dots,d-1\}$, with the sign $\epsilon=1$ if $x_{n/d}=x=0$ and
  $\epsilon=-1$ if $x_{n/d}=x=1$.
\end{lemma}
\begin{proof}
Let us prove that the described set $\nuke$ of elements of $F_n$
is the nucleus of the recursion $\Psi$.

It follows directly from the wreath recursion that if
$\Psi(g)=\pair<g_0,g_1>\sigma^{k}$, then
$\|g_0\|+\|g_1\|\le\|g\|$. Let us call an element $g$
\emph{stable} if $\|g_i\|=\|g\|$ for some $i\in\{0, 1\}$ and the
element $g_i$ is also stable. In other words, an element $g\in
F_n$ is stable if for every $n$ there exists a word $v_n\in X^n$
(which is necessarily unique) such that $\|g|_{v_n}\|=\|g\|$. Note
that then for every $u\in X^n$ different from $v_n$ we have
$g|_u=1$, and $v_n$ is the prefix of $v_{n+1}$. Hence there exists
an infinite sequence $w_g$ such that $v_n$ is the beginning of
length $n$ of $w_g$. We call the sequence $w_g$ the \emph{kneading
sequence} of $g$.

If we prove that the set of stable elements is finite, we will prove
that the recursion is contracting and the nucleus is the set of stable
elements, since the length of elements is non-increasing under taking
restrictions.

It is easy to see that if $g$ is stable, then $g^{-1}$ is also
stable and $w_{g^{-1}}=(w_g)^g$. If $g$ is stable and
$g=a_{i_1}^{\epsilon_1}a_{i_2}^{\epsilon_2}\cdots
a_{i_k}^{\epsilon_k}$ for $\epsilon_i\in\{\pm 1\}$ and $k=|g|$,
then $w_g$ is also the kneading sequence of
$a_{i_1}^{\epsilon_1}$, $(w_g)^{a_{i_1}^{\epsilon_1}}$ is the
kneading sequence of $a_{i_2}^{\epsilon_2}$,
$(w_g)^{a_{i_1}^{\epsilon_1}a_{i_2}^{\epsilon_2}}$ is the kneading
sequence of $a_{i_3}^{\epsilon_3}$, etc. Otherwise the length of
$g$ will decrease along $w_g$. Hence the set of kneading sequences
of the stable elements coincides with the set $K$ of the kneading
sequences of the generators and their inverses. Let $\Gamma$ be
the graph with set of vertices $K$ in which for every $a_i$
there is an edge from the kneading sequence $w_{a_i}$ to the
kneading sequence $w_{a_i^{-1}}=(w_{a_i})^{a_i}$ labeled by $a_i$
(if we go in the opposite direction then the edge is labeled by
$a_i^{-1}$). Then for every stable element
$g=a_{i_1}^{\epsilon_1}a_{i_2}^{\epsilon_2}\cdots
a_{i_k}^{\epsilon_k}$ there exists a path (without back-tracking)
in the graph $\Gamma$ labeled by $a_{i_1}^{\epsilon_1},
a_{i_2}^{\epsilon_2}, \ldots, a_{i_k}^{\epsilon_k}$.

It follows from the definitions that the kneading sequence of
$a_i$ is $(x_ix_{i-1}\ldots x_1 1 x_{n-1}\ldots x_{i+1})^\omega$
and the kneading sequence of $a_i^{-1}$ is $(x_ix_{i-1}\ldots x_1
0 x_{n-1}\ldots x_{i+1})^\omega$. Let, as above, $d$ be the
maximal number such that $vx=u^d$ for some $x\in\{0, 1\}$. We
denote by $x'\ne x$ the other letter of the alphabet. If $x=1$,
then the kneading sequence of $a_i$ is of period $n/d$ and the
kneading sequence of $a_i^{-1}$ is of period $n$. If $x=0$, then
the period of the kneading sequence of $a_i$ is $n$ and the period
of the kneading sequence of $a_i^{-1}$ is $n/d$. In each case the
set $K$ of kneading sequences of stable elements of $\kn{v}$
contains $n+n/d$ elements. The graph $\Gamma$ consists of $n/d$
disjoint ``stars'': every sequence of period $n$ is connected to
the sequence of period $n/d$ obtained by changing the respective
$x'$ to $x$. Every sequence of period $n$ is hence a vertex of
$\Gamma$ of valence $1$ and every sequence of period $n/d$ is a
vertex of valence $d$. Consequently, paths without back-tracking
have length at most 2 and the labels read on these paths are the
elements described in the lemma. The rest of the lemma easily
follows.
\end{proof}

\begin{proposition}\label{prop:kv:abelian}
  Let $e_1, \ldots e_n$ be a free basis of the group $\Z^n$.  The map
  $\rho:F_n\to\Z^n$ given by $a_i\mapsto e_i$ can be extended to a
  homomorphism $\rho:\kv\to\Z^n$.  Consequently,
  $\kv/\kv'\cong\Z^n$.
\end{proposition}
\begin{proof}
  Write $\pi:F_n\to\kv$ the natural quotient map.  Assume by
  contradiction that $\rho$ does not factor through $\pi$; then there
  exists $g\in\ker\pi\setminus\ker\rho$. Let $g=a_{i_1}^{\pm
    1}a_{i_2}^{\pm 1}\cdots a_{i_k}^{\pm 1}$ be a shortest (i.e., one
  with smallest $k$) element of $F_n$ in $\ker\pi\setminus\ker\rho$.
  Clearly we may assume $k\ge2$, since every generator has an active
  state, and is therefore non-trivial in $\autxs/\autxs'$.

  Let us write $\Psi(g)=\pair<g_0,g_1>$. Then, by considering the
  formulae, we see that the sum of the lengths of $g_0$ and $g_1$ is
  at most equal to $k$; it could be shorter, if some cancellation
  inside $g_0$ or $g_1$ occurs. The formulae~\eqref{eq:kv:Psi} also
  show that $\rho(g_0)+\rho(g_1)=T(\rho(g))$, where $T:\Z^n\to\Z^n$ is
  the isomorphism $e_i\mapsto e_{i-1\mod n}$. Consequently, either
  $\rho(g_0)$ or $\rho(g_1)$ is non-trivial. Since $\kv$ is
  contracting, $g_0$ and $g_1$ are shorter than $g$, and we have
  contradicted the minimality of $g$.
\end{proof}

\begin{lemma}\label{lem:kv:2}
  Define the subgroups $\eker_i\le F_n$ inductively by $\mathcal
  E_0=1$ and
  \[\eker_{i+1}=\Psi^{-1}(\eker_i\times\eker_i).
  \]
  Then $\eker_\infty=\bigcup_{i=1}^\infty\eker_i$ is the
  kernel of the natural homomorphism $\pi:F_n\arr\kv$.
\end{lemma}
\begin{proof}
  By contraction, for every $g\in F_n$ there exists $k$ such that all
  coordinates of $\Psi^k(g)$ belong to the nucleus. The intersection
  of the nucleus with $\ker(\pi)$ is trivial.  Therefore, if $g$
  belongs to $\ker(\pi)$, then there exists $k$ such that all the
  coordinates of $\Psi^k(g)$ are trivial, i.e.,
  $g\in\ker\Psi^k=\eker_k$.
\end{proof}

\subsection{Endomorphisms}
Let us define the following endomorphism $\varphi$ of the free
group $F_n$:
\begin{equation}\label{eq:kv:phi}
  \varphi(a_n)=a_1^2, \qquad\varphi(a_i)=\left.\begin{cases}
    a_{i+1} & \text{ if }x_i=0\\
    a_1a_{i+1}a_1^{-1} & \text{ if }x_i=1\end{cases}\right\}\text{ when }1\le i<n.
\end{equation}

We will also write $\varphi(g)$ for its image in $\kv$. It follows
directly from the definition that
\[
\varphi(a_i)=\pair<a_i, 1>\text{ when }1\le i<n,\qquad\varphi(a_n)=\pair<a_n, a_n>.
\]

\begin{proposition}\label{prop:kv:varphi}
  The substitution $\varphi$ induces an endomorphism of the group
  $\kv$ such that
  \[
  \varphi(g)=\pair<g, \rho_n(g)>,
  \]
  where $\rho_n$ is the composition of the abelianization $a_i\mapsto
  e_i\in\Z^n$ and the projection $e_n\mapsto a_n$, $e_i\mapsto 1$ for
  $i\ne n$.
\end{proposition}

\begin{lemma}\label{lem:kv:3}
  For every $\ell\ge 1$ we have
  \[
  \eker_{\ell+1}\le\eker_1\cdot\varphi(\eker_\ell)\cdot
  a_1^{-1}\varphi(\eker_\ell)a_1.
  \]
\end{lemma}
\begin{proof}
  Consider $g\in\eker_{\ell+1}$. Then $\Psi(g)=\pair<g_0, g_1>$, with
  $g_0, g_1\in\eker_\ell$. We have $\Psi\varphi(g_i)=\pair<g_i,
  \rho_n(g_i)>$ for all $i\in\{0,1\}$.  Moreover, $\rho_n(g_i)=1$,
  since $g_i=1$ in $\kv$ and the image of $\rho$ is free. Hence,
  \[
  \Psi(g)=\pair<g_0, g_1>=\Psi(\varphi(g_0)\cdot a_1^{-1}\varphi(g_1)a_1),
  \]
  so $g=h\cdot\varphi(g_0)\cdot a_1^{-1}\varphi(g_1)a_1$ for some
  $h\in\eker_1$.
\end{proof}

\subsection{L-presentation}
We give a recursive presentation of $\kv$ by generators and relations.

For $r\ge0$ define the set of commutators
\begin{multline}
  \label{eq:kv:R}
  \rel_r=\bigg\{\Big[a_i,a_j^{a_1^{2k}}\Big]:\,2\le i,j\le n,\,0\le2k\le
  r,\,x_{i-1}\neq x_{j-1}\bigg\}\\
  \cup\bigg\{\Big[a_i,a_j^{a_1^{2k+1}}\Big]:\,2\le i,j\le n,\,0\le2k+1\le
  r,\,x_{i-1}=x_{j-1}\bigg\}.
\end{multline}

\begin{lemma}\label{lem:kv:4}
  The subgroup $\eker_1< F_n$ is the normal closure of $\mathcal
  R_\infty=\bigcup_{r\ge0}\rel_r$.
\end{lemma}
\begin{proof}
  It is easy to see that $\rel_r\subset\eker_1$ for all
  $r$. Let now $g\in\eker_1$ be any non-trivial reduced group
  word. We can write it in the form $a_1^kc_1c_2\ldots c_m$, where
  $c_s=a_{i_s}^{\pm a_1^{j_s}}$ for some $i_s\in\{2,\dots,n\}$ and
  some $j_s\in\Z$. Since $g=1$ in $\kv$, we have $k=0$ by
  Proposition~\ref{prop:kv:abelian}. Let us write
  $\Psi(g)=\pair<g_0,g_1>$. Then each of $g_0,g_1$ respectively is the
  product of some $a_{i_s-1}^{\pm a_n^{k_s}}$, where $|2k_s-j_s|\le1$
  and only those $s$ with $j_s$ congruent modulo 2, respectively not
  congruent, to $x_{i_s-1}$, are selected.

  Now since $g_0=g_1=1$, there must exist $s<t$ such that $i_s=i_t$
  and $k_s=k_t$ and $c_s,c_t$ have opposing signs, both of them
  occurring in the same $g_i$, and none of the $c_{s+1},\dots,c_{t-1}$
  contributing to that $g_i$. Then the relations in $\mathcal
  R_\infty$ allow the commutation of $c_s$ with
  $c_{s+1},\dots,c_{t-1}$, and its eventual cancellation with
  $c_t$. One then proceeds by induction on the length of $g$.
\end{proof}

\begin{theorem}
  The group $\kv$ has the following presentation:
  \[
  \kv=\Big\langle a_1,\ldots a_n\,\Big|\,\varphi^\ell(\rel_2)\text{
    for all }\ell\ge0\Big\rangle,
  \]
  where $\varphi$ and $\rel_2$ and given respectively
  in~\eqref{eq:kv:phi} and~\eqref{eq:kv:R}.
\end{theorem}
\begin{proof}
  In view of Lemmata~\ref{lem:kv:2}, \ref{lem:kv:3}
  and~\ref{lem:kv:4}, it suffices to show for all $r\ge0$ that
  $\rel_r$ is a consequence of $\mathcal
  R^*=\bigcup_{k\ge0}\varphi^k(\rel_2)$. Consider therefore a
  relation $w_{ijk}=[a_i,a_j^{a_1^k}]$ with $k\ge3$. We write $g\equiv
  h$ to mean that they are equivalent under relations in $\mathcal
  R^*$.

  We first write $a_j^{a_1^k}=a_j^{a_1^{k-2}}[a_j,a_1^2]^{a_1^{k-2}}$;
  then by $[a,bc]=[a,c][a,b]^c$ the relation $w_{ijk}$ follows from
  the relations $w_{ij(k-2)}$ and $[a_i,[a_j,a_1^2]^{a_1^{k-2}}]$.

  We then note that $[a_p,[a_p,a_1]]\in\rel_2$ for all $p$, and more
  generally $[a_p,[a_p,a_q]]\in\rel^*$ for all $p\ge q$: if
  $x_{p-1}=x_{q-1}$, this relation is obtained as
  $\varphi([a_{p-1},[a_{p-1},a_{q-1}]])$, while if $x_{p-1}\neq
  x_{q-1}$, it follows from $[a_p,a_q]\in\rel_2$. In particular
  $[a_n,[a_n,a_{i-1}]]\in\rel^*$, and $[a_1^2,[a_1^2,a_i]]\in\rel^*$.
  This allows us to reduce $[a_i,[a_j,a_1^2]^{a_1^{k-2}}]$ to
  $[a_i,[a_j,a_1^2]^e]$ for some $e\in\{1,a_1^{-1}\}$.  Now if $e=1$
  then we have $[a_i,a_j,a_1^2]\equiv w_{ij2}[a_i,a_j]$, and if
  $e=a_1^{-1}$ then we have
  $[a_i,[a_j,a_1^2]^{a_1^{-1}}]=w_{ij1}^{-a_1^{-1}}w_{ij1}^{-a_1^{-1+2a_j-1}}$.
\end{proof}

Since the endomorphism $\varphi$ is injective, we can embed the group
$\kv$ into its ascending HNN-extension by $\varphi$, i.e.\ in the
group generated by $\kv$ and an element $t$ whose action by
conjugation on $\kv$ coincides with $\varphi$.

Setting $a=a_1^{-1}$, we have a new generating system $\{a, a^t,
a^{t^2}, \ldots, a^{t^{n-1}}\}$ of the group $\kv$; the
identification is
\begin{equation}\label{eq:kv:identify}
  \begin{aligned}
    a_{i+1} &= a^{x_i}a_i^ta^{-x_i} = a^{x_i}a^{x_{i-1}t}\cdots
    a^{x_1t^{i-1}}a^{-t^i}a^{-x_1t^{i-1}}\cdots a^{-x_{i-1}t}a^{-x_i},\\
    a^{t^i} &= a_{i+1}^{-a^{x_i+x_{i-1}t+\dots+x_1t^{i-1}}}.
  \end{aligned}
\end{equation}

The last theorem then yields
\begin{theorem}
  Write $p(t)=x_{n-1}t+x_{n-2}t^2+\dots+x_1t^{n-1}\in\Z[t]$. Then the
  group $\kv$ is isomorphic to the subgroup $\langle a, a^t, a^{t^2},
  \ldots, a^{t^{n-1}}\rangle$ of the finitely presented group
  \[
  \Big\langle a, t\,\Big|\, a^{t^n-2a^{p(t)}},
  \big[a^{t^i},a^{t^ja}\big], \big[a^{t^i},a^{t^ja^3}\big]\text{ for
    all }1\le i,k<n\Big\rangle.
  \]
\end{theorem}
\begin{proof}
  The first relation is $a_n^t=a_1^2$. The others are obtained by
  rewriting $[a_i,a_j^{a_1^k}]$ for $0\le k\le 2$ in terms of $a$ and
  $t$. Indeed if $x_i=x_j$ then the relation $[a_{i+1},a_{j+1}^{a_1}]$
  yields
  \[[a^{x_{i-1}t}\cdots a^{x_1t^{i-1}}a^{-t^i}a^{-x_1t^{i-1}}\cdots a^{-x_{i-1}t},(a^{x_{j-1}t}\cdots a^{x_1t^{j-1}}a^{-t^j}a^{-x_1t^{j-1}}\cdots a^{-x_{j-1}t})^{a_1}],\]
  which is equivalent to $[a^{t^i},a^{t^ja}]$ for all $i,j$. If
  $x_i\neq x_j$, then the relation $[a_{i+1},a_{j+1}]$ gives the same
  commutation relation as above, while $[a_{i+1},a_{j+1}^{a_1^2}]$
  gives the equivalent form $[a^{t^i},a^{t^ja^3}]$.
\end{proof}

\subsection{Torsion}
We show that $\kv$ is torsion-free for all $v$:

\begin{proposition}
  Let $\rho:\kv\to\kv/\kv'\cong\Z^n$ be the canonical epimorphism.
  An element $g\in\kv$ is level-transitive if and only if all
  coordinates of $\rho(g)$ are odd.
\end{proposition}
\begin{proof}
  We write $\rho(a_i)=e_i$, where $(e_i)$ is a free basis of $\Z^n$.

  Recall that $\autxs/\autxs'$ is $(\ZZ2)^\omega$, the identification
  being $\tau:g\mapsto(i_0,i_1,\dots)\in(\ZZ2)^\omega$ where $i_m$ is
  the parity of the number of active restrictions $g|_v$ for all $v\in
  X^m$. A tree automorphism is \emph{active} if it acts non-trivially
  on the first level.

  It is well known that $g\in\autxs$ is level-transitive if and only
  if $\tau(g)=(1,1,\dots)$. Now if $g\in\kv$, then $\tau(g)_m$ is the
  parity of the exponent of $e_{m+1\bmod n}$ in $\rho(g)$.
\end{proof}

\begin{proposition}\label{prop:kv:torsionfree}
  The group $\kv$ is torsion free.
\end{proposition}
\begin{proof}
  Apply Proposition~\ref{prop:torsionfree} with $N=G'$.
\end{proof}

\subsection{Weak branchness}
We show that the groups $\kv$ are weakly branch and residually poly-$\Z$:

\begin{theorem}
  Let $v\in\xs$ be a non-empty sequence. Then the group $\kv$ is
  weakly branch on $\kv'$.
\end{theorem}
Note that the group $\kn{\emptyset}$ is isomorphic to $\Z$ and so is
not weakly branch.
\begin{proof}
  As a first step, let us check that $\kv'\times\kv'\le\kv$. It is
  easy to check that
  \[\left\langle a_1, \ldots, a_{n-1}, a_1^{a_n},
    \ldots, a_{n-1}^{a_n}\right\rangle\times\left\langle a_1, \ldots,
    a_{n-1}, a_1^{a_n}, \ldots, a_{n-1}^{a_n}\right\rangle\le\kv.
  \]
  But
  \[
  \kv'\le\left\langle a_1, \ldots, a_{n-1}, a_1^{a_n}, \ldots,
    a_{n-1}^{a_n}\right\rangle,
  \]
  since $\kv'$ is contained in the kernel of the homomorphism
  $\rho_n$. Hence $\kv'\times\kv'\le\kv$.

  Let now $g=\pair<g_0, g_1>$ be an arbitrary element of the subgroup
  $\kv'\times\kv'$. Then $\varphi(g_0)=\pair<g_0,
  \rho_n(g_0)>\in\kv'$ and
  $a_1^{-1}\varphi(g_1)a_1=\pair<\rho_n(g_1), g_1>\in\kv'$, where
  $\varphi$ and $\rho_n$ are as in Proposition~\ref{prop:kv:varphi}. But
  $\rho_n(g_0)=\rho_n(g_1)=1$, since $g_0, g_1\in\kv'$. Therefore,
  $g=\varphi(g_0)\cdot a_1^{-1}\varphi(g_1)a_1\in\kv'$.

  It suffices now to check that $\kv'$ is not trivial. For that
  purpose, note that $[a_1,a_2]=\pair<a_1^{\pm a_n},a_1^{\mp1}>$ has
  non-trivial image in $\Z^n$.
\end{proof}

\begin{proposition}
  We have $\kv'/\kv'\times\kv'\cong\Z^{n-1}$.
\end{proposition}
\begin{proof}
  Since the $[a_i,a_j]$ belong to $\kv'\times\kv'$ if $i,j\ge2$, we
  see that $\kv'/\kv'\times\kv'$ is abelian and generated by the
  $[a_1,a_{i+1}]$ for $i\in\{1,\dots,n-1\}$.

  We then have $[a_1,a_{i+1}]=\pair<a_i^{\mp a_n},a_i^{\pm 1}>$; by
  Proposition~\ref{prop:kv:abelian}, these are independent in
  $(\kv/\kv')^2$.
\end{proof}

It follows that the groups $\kv$ are not branch; otherwise, they would
be branch on a subgroup $K$ with $\kv'\le K\le\kv$, the last inclusion
of finite index. Then $K/\kv'$ has rank $n$, so $K/(\kv'\times\kv')$
has Hirsch length $2n-1$, but it contains
$(K\times K)/(\kv'\times\kv')$ of Hirsch length $2n$, a contradiction.

It follows that $\kv$ admits arbitrarily large poly-$\Z$ quotients:
the successive quotients along the descending series
$\kv>\kv'>\kv'\times\kv'>\dots$ are all free abelian.

\begin{corollary}
  The group $\kv$ is left orderable.
\end{corollary}

\section{The groups $\kwv$}\label{ss:kwv}
Let $w=y_1\ldots y_k\in\xs$ and $v=x_1\ldots x_n\in\xs$ be a pair of
non-empty words such that $y_k\ne x_n$. We denote by $\kwv$ the
subgroup of $\autxs$ generated by the elements $b_1, \ldots b_k, a_1,
\ldots, a_n$ defined by
\begin{gather*}
  b_1=\sigma,\qquad
  b_{j+1}=\left.\begin{cases}\pair<b_j, 1> & \text{ if }y_j=0\\
    \pair<1, b_j> & \text{ if }y_j=1\end{cases}\right\}\text{ when }1\le j<k,\\
  a_1=\begin{cases}\pair<b_k, a_n> & \text{ if }y_k=0\text{ and }x_n=1,\\
    \pair<a_n, b_k> & \text{ if }y_k=1\text{ and }x_n=0,\end{cases}\qquad
  a_{i+1}=\left.\begin{cases}\pair<a_i, 1> & \text{ if }x_i=0\\
    \pair<1, a_i> & \text{ if }x_i=1\end{cases}\right\}\text{ when }1\le i<n.
\end{gather*}

\begin{figure}
  \begin{center}
    \includegraphics{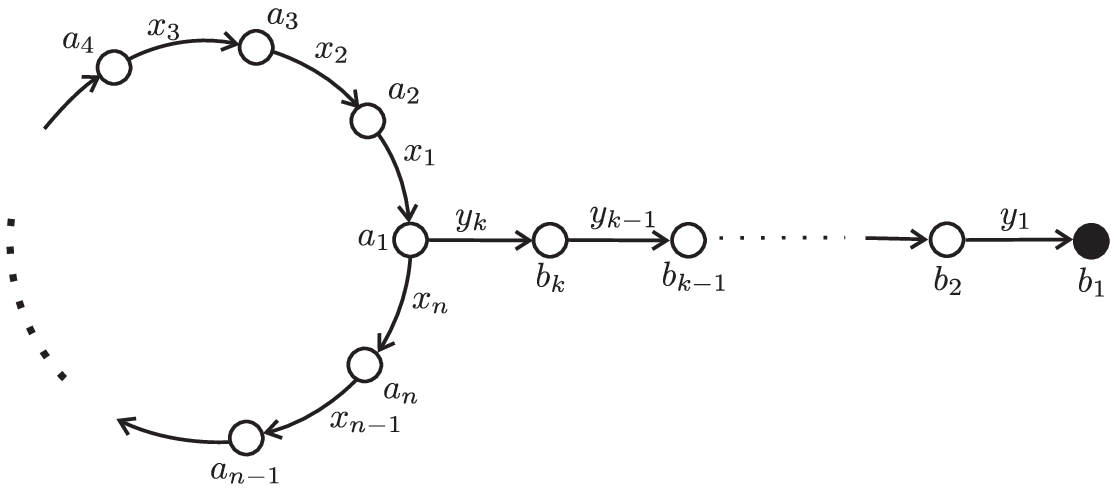}
    \caption{Automaton generating $\kwv$}\label{fig:kuv}
  \end{center}
\end{figure}

In other words, the group $\kwv$ is generated by the automaton whose
Moore diagram is shown on Figure~\ref{fig:kuv}.  In this diagram, only
the edges leading to non-trivial states are drawn, and the active state
$b_1$ is labelled by $\sigma$.

Let us denote by $v'$ the word obtained from the word $v$ by changing
$0$ to $1$ and $1$ to $0$.  Recall from~\S\ref{ss:kv} that the
automorphism $\alpha\in\autxs$ defined by the recursion
$\alpha=\pair<\alpha,\alpha>\sigma$ exchanges $v$ and $v'$, and
conjugates $\kv$ into $\kn{v'}$.  Similarly, the groups $\kwv$ and
$\kn{w',v'}$ are conjugate under $\alpha$.

By convention, let us accept that $b_{k+1}=a_{n+1}=a_1$.

\begin{lemma}
  The group $\kwv$ is recurrent and level-transitive.
\end{lemma}
\begin{proof}
  The projections of the $a_{i+1}$, $a_{i+1}^{b_1}$, $b_{j+1}$ and
  $b_{j+1}^{b_1}$, all fixing the first level of $\xs$, give all
  generators; the action is therefore recurrent.

  To check that the action of a recurrent group is level-transitive,
  it suffices to check that it is transitive on the first level; this
  is achieved by $b_1$.
\end{proof}

\subsection{Wreath recursion}
We aim to give a recursive presentation of $\kwv$ by generators
and relations. Until the end of this section, we assume that
$(k,n)\neq(1,1)$; indeed if $k=n=1$ then $\kwv$ is infinite
dihedral, and has to be taken care of separately.

\begin{proposition}\label{prop:kwv:abelian}
  The generators $a_i,b_j$ of $\kwv$ have order $2$, and its
  abelianization $\kwv/\kwv'$ is $(\ZZ2)^{k+n}$, generated by the images of
  the $a_i$ and $b_j$.
\end{proposition}
\begin{proof}
  It is clear that the generators are involutions. The image of $\kwv$
  in $\autxs/\autxs'$ is generated by the infinite pre-periodic vectors
  $(0,\dots,1,\overline 0)$ corresponding to the images of $b_j$, and
  by $(0,\dots,0,\overline{0,\dots,1,\dots,0})$ corresponding to the
  images of $a_i$; these $k+n$ vectors are linearly independent.
\end{proof}

Let $d$ be maximal such that $v=u^d$ for some $u\in\xs$. We have
then, considering the indices modulo $n$,
\[\Psi([a_i, a_{i+j(n/d)}])=\pair<1,[a_{i-1},a_{i-1+j(n/d)}]>\text{ or }
  \pair<[a_{i-1},a_{i-1+j(n/d)}],1>
\]
depending on whether $x_{i-1}=x_{i-1+j(n/d)}$ or not.  It follows that
$a_i$ and $a_{i+j(n/d)}$ commute for every $i\in\{1,\dots,n\}$ and
every $j\in\{1,\dots,d-1\}$.

Let $F$ be the group given by the presentation
\[F=\langle a_1, \ldots, a_n, b_1, \ldots, b_k\;|\;a_i^2, b_j^2,
[a_i,a_{i+\ell(n/d)}]\text{ for all }1\le i\le n,\;1\le j\le k,\;1\le\ell<d\rangle.
\]
It is the free product of $k$ copies of $\ZZ2$ and $n/d$ copies of
$(\ZZ2)^d$.  Let us denote the subgroup $\langle a_i, a_{i+n/d},
\ldots, a_{i+(d-1)n/d}\rangle\cong(\ZZ2)^d$ by $A_i$.

Let $\Psi:F\arr F\wr\symm$ be given by the recursive definition of
the group $\kwv$, i.e.,
\begin{equation}\label{eq:kwv:Psi}
  \begin{aligned}
    \Psi(b_1)=\sigma,\qquad \Psi(b_{j+1}) =
    \left.\begin{cases}\pair<b_j, 1> & \text{ if }y_j=0\\ \pair<1,
        b_j> & \text{ if }y_j=1\end{cases}\right\}\text{ when }1\le j<k,\\
    \Psi(a_1)=\begin{cases}\pair<b_k, a_n> & \text{ if }x_n=1,\\
      \pair<a_n, b_k> & \text{ if }x_n=0,\end{cases}\qquad
    \Psi(a_{i+1})=\left.\begin{cases}\pair<a_i, 1> & \text{ if }x_i=0\\
      \pair<1, a_i> & \text{ if }x_i=1\end{cases}\right\}\text{ when }1\le i<n.
  \end{aligned}
\end{equation}

We can iterate the map $\Psi$ to obtain a homomorphism
$\Psi^k:F\arr F\wr\symm\wr\cdots\wr\symm$, where $\symm$
appears $k$ times in the iterated permutational wreath product.

\begin{lemma}\label{lem:kwv:1}
  The wreath recursion $\Psi$ is contracting. Its nucleus is
  $\{b_1,\dots,b_k\}\cup\bigcup A_i$.
\end{lemma}

\begin{proof}
Let us denote by $\|g\|$ the minimal number of generators $a_i$ in
a representation of $g$ as a product of generators.

It follows from the definition of the recursion $\Psi$ that for
every $g=\pair<g_0,g_1>\sigma^k\in\kwv$ we have $\|g\|\ge
\|g_0\|+\|g_1\|$. As in the proof of Lemma~\ref{lem:kv:1}, we say
that $g\in\kwv$ is \emph{stable} if there exists an infinite
sequence $w_g$ (the \emph{kneading sequence} of $g$) such that for
every $n$ we have $\|g\|=\|g|_{v_n}\|$ where $v_n$ is the
beginning of length $n$ of $w_g$. It also follows that for any
word $v\in X^n$ different from $v_n$ the restriction $g|_v$
belongs to $\langle b_i\rangle_{i=1, \ldots, k}$, and hence that there
exists $m\in\N$ such that $g|_{vw}=1$ for all $w\in X^m$.

Consequently, if we prove that the set of stable elements is finite,
then we show that the recursion is contracting and the nucleus is the
set of stable elements and their restrictions.

As in the proof of Lemma~\ref{lem:kv:1}, consider the set $K$
of kneading sequences of the (stable) generators $a_i$, and the
corresponding graph $\Gamma$. Here also the labels along the paths of
the graph $\Gamma$ correspond to the stable elements.

The kneading sequence of $a_i$ is equal to
$w_{a_i}=(x_{i-1}x_{i-2}\ldots x_1x_nx_{n-1}\ldots x_{i})^\omega$. If
$d$ is as above, then the kneading sequence of $a_i$ is of period
$n/d$ and the graph $\Gamma$ consists of $n/d$ vertices with $d$ loops
attached to every vertex. The labels of the loops attached to vertex
$w_{a_i}$ are $a_i, a_{i+n/d}, \ldots, a_{i+(d-1)n/d}$.  Consequently,
every stable element belongs to the finite group $A_i$ for some
$i\in\{1,\dots,n/d\}$. The rest of the lemma follows.
\end{proof}

\begin{lemma}\label{lem:kwv:2}
  Define the subgroups $\eker_\ell\le F$ inductively by $\mathcal
  E_0=1$ and
  \[\eker_{\ell+1}=\Psi^{-1}(\eker_\ell\times\eker_\ell).
  \]
  Then $\eker_\infty=\bigcup_{\ell=1}^\infty\eker_\ell$ is the
  kernel of the natural homomorphism $F\arr\kwv$.
\end{lemma}
\begin{proof}
  Similar to Lemma~\ref{lem:kv:2}.
\end{proof}

\subsection{Endomorphisms}
In this subsection we define, for each $\kwv$, an endomorphism
$\varphi$ of $F$, and elements $s,t\in F$ and $m\in\{1,2\}$, with the
following meaning: $\Psi(a_1^t)=\pair<b_k,a_n^s>$; in $\kwv$ the
subgroup $\langle a_n^s,b_k\rangle$ is dihedral of order $2^{m+1}$;
and the subgroup $\langle a_1^t,b_1\rangle$ is dihedral of order
$2^{m+2}$; and $\varphi(s)\in\{t,b_1t\}$. We use these to construct an
endomorphism of $\kwv$, considering different cases in turn and taking
care of small values of $k$ and $n$.

Assume first that $k,n\ge2$. Then $[a_n^s,b_k]$ is trivial in $\kwv$
for at least one value of $s\in\langle b_1\rangle$. Let $t\in\langle
b_1,b_2\rangle$ be such that $\Psi(a_1^t)=\pair<b_k,a_n^s>$. Set $m=1$.

Assume next that $k\ge3$ and $n=1$.  Then there exists $r\in\langle
b_1\rangle$ such that $[b_k^r,b_{k-1}]=1$ in $\kwv$; so there exists
$s\in\langle b_1,b_2\rangle$ such that $[a_1^s,b_k]=1$ in $\kwv$; so
there exists $t\in\langle b_1,b_2,b_3\rangle$ such that
$\Psi(a_1^t)=\pair<b_k,a_1^s>$. Set also $m=1$.

Assume next that $k=2$ and $n=1$. Then there exists $r,s\in\langle
b_1\rangle$ such that $b_2^r=\pair<b_1,1>$ in $\kwv$ and
$a_1^{rs}=\pair<b_2,a_1>$ in $\kwv$; so there exists $t\in\langle
b_1,b_2\rangle$ such that $a_1^t=\pair<b_2,a_1^s>$. Set $m=2$.

Finally, consider the case $k=1$ and $n\ge2$. Set $s=1$, and let
$t\in\langle b_1\rangle$ be such that $a_1^t=\pair<b_1,a_n>$. Set $m=2$.

Now define the endomorphism $\varphi$ of $F$ by
\begin{gather}\label{eq:kwv:phi}
  \begin{aligned}
    \varphi(b_j)=\left.\begin{cases} b_{j+1} & \text{ if }y_j=0\\
      b_{j+1}^{b_1} & \text{ if }y_j=1\end{cases}\right\}\text{ when }1\le
    j<k,&\qquad \varphi(b_k)=a_1^t,\\
    \varphi(a_i)=\left.\begin{cases} a_{i+1} & \text{ if }x_i=0\\
      a_{i+1}^{b_1} & \text{ if }x_i=1\end{cases}\right\}\text{ when }1\le
    i<n,&\qquad \varphi(a_n)=\begin{cases} a_1 & \text{ if }x_n=0,\\
      a_1^{b_1} & \text{ if }x_n=1.\end{cases}
  \end{aligned}\\
\intertext{We then have}
  \begin{aligned}
    \Psi\varphi(b_j)=\pair<b_j,1>\text{ when }1\le j<k,
    &\qquad\Psi\varphi(b_k)=\pair<b_k,a_n^s>,\\
    \Psi\varphi(a_i)=\pair<a_i,1>\text{ when }1\le i<n,
    &\qquad\Psi\varphi(a_n)=\pair<a_n,b_k>.
  \end{aligned}
\end{gather}

\begin{proposition}\label{prop:kwv:varphi}
  The endomorphism $\varphi:F\to F$ induces an endomorphism of the
  group $\kwv$ such that
  \[
  \varphi(g)=\pair<g, \rho(g)>,
  \]
  where $\rho$ is an endomorphism with finite image; more
  precisely, $\rho(b_j)=\rho(a_i)=1$ for $j<k$ and $i<n$; and
  $\rho(b_k)=a_n^s$ and $\rho(a_n)=b_k$, so that $\rho(G)$ is dihedral
  of order $2^{m+1}$.
\end{proposition}
\begin{proof}
  It is clear that $\rho$ is an endomorphism if $m=1$, since then it
  factors through $\kwv\to\autxs/\autxs'$.

  Consider then the case $k=1$. The map $\rho$ then factors through
  $\kwv\to\kwv\wr\symm\to(\autxs/\autxs')\wr\symm$.

  Consider finally the case $k=2,n=1$. The map $\rho$ can then be
  seen, by direct calculation, to factor through $\kwv\to\aut(X^4)$.
\end{proof}

Define now $\overline F=F/\langle
\varphi^\ell(a_n^sb_k)^{2^m}:\,\ell\in\N\rangle$, and denote by
$\overline{\eker_\ell}$ the image of $\eker_\ell$ in $\overline F$.
Again $\varphi$ induces an endomorphism of $\overline F$.

\begin{lemma}\label{lem:kwv:3}
  For every $\ell\ge 1$ we have
  \[
  \overline{\eker_{\ell+1}}\le\overline{\eker_1}\cdot\varphi(\overline{\eker_\ell})\cdot
  b_1\varphi(\overline{\eker_i})b_1.
  \]
\end{lemma}
\begin{proof}
  The proof is similar to that of Lemma~\ref{lem:kv:3}.  Consider
  $g\in\overline{\eker_{\ell+1}}$. Then $\Psi(g)=\pair<g_0, g_1>$,
  with $g_0, g_1\in\overline{\eker_\ell}$. We have
  $\Psi\varphi(g_i)=\pair<g_i, \rho(g_i)>$ for all $i\in\{0,1\}$,
  where $\rho$ is the endomorphism with finite image given by
  Proposition~\ref{prop:kwv:varphi}. Moreover, $\rho(g_i)=1$, since
  $\rho(g_i)$ is a relation in $\kwv$. Hence,
  \[
  \Psi(g)=\pair<g_0, g_1>=\Psi(\varphi(g_0)\cdot b_1\varphi(g_1)b_1),
  \]
  and $g=h\cdot\varphi(g_0)\cdot b_1\varphi(g_1)b_1$ for some
  $h\in\overline{\eker_1}$.
\end{proof}

\subsection{L-presentation}
Let $O$ denote the elements of $\langle a_1^t,b_1\rangle<\overline F$
that contain an odd number of $b_1$'s, and let $E$ denote those
elements that contain an even number of $b_1$'s. Both of these sets
are finite. Define then
\begin{align}
  \label{eq:kwv:R}
  \rel&=\Big\{[b_i,b_j^w]:\,2\le i,j\le k,\,w\in
  O\text{ if }y_{i-1}=y_{j-1},\,w\in E\text{ if }y_{i-1}\neq
  y_{j-1}\Big\}\notag\\
  &\cup\Big\{[a_i,b_j^w]:\,2\le i\le n,\,2\le j\le k,\,w\in O\text{ if
  }x_{i-1}=y_{j-1},\,w\in E\text{ if }x_{i-1}\neq y_{j-1}\Big\}\\
  &\cup\Big\{[a_i,a_j^w]:\,2\le i,j\le n,\,w\in
  O\text{ if }x_{i-1}=x_{j-1},\,w\in E\text{ if }x_{i-1}\neq
  x_{j-1}\Big\}.\notag
\end{align}

\begin{lemma}\label{lem:kwv:4}
  We have $\overline{\eker_1}\le \rel^{\overline
    F}\le\overline{\eker_\infty}\le \overline F$.
\end{lemma}
\begin{proof}
  It is easy to see that $\rel\subset\overline{\eker_\infty}$. Let now
  $g\in\overline{\eker_1}$ be any non-trivial reduced group word. We
  can write it in the form $hc_1c_2\ldots c_m$, where
  $c_s=a_{i_s}^{h_s}$ or $=b_{j_s}^{h_s}$ for some $i_s,j_2\ge2$ and
  some $h,h_s\in\langle a_1^t,b_1\rangle$.  Since $g=1$ in $\kv$, we
  have $h=1$; this is clear if $a_1^t$ and $b_1$ commute, because then
  $\langle a_1^t,b_1\rangle$ embeds in $\kwv/\kwv'$, by
  Proposition~\ref{prop:kwv:abelian}. On the other hand, if
  $h=(a_1^tb_1)^2\neq1$, then the first projection of $g$ contains
  precisely two $a_n$'s that do not cancel, and there will exist two
  vertices such that one of the $a_n$'s projects on each; then that
  state will be non-trivial, again using
  Proposition~\ref{prop:kwv:abelian}, so $g\neq1$, a contradiction.

  Let us write $\Psi(g)=\pair<g_0,g_1>$. Then each of $g_0,g_1$
  respectively is the product of some $a_{i_s-1}^{k_s}$ and
  $b_{j_s-1}^{k_s}$, where $k_s\in\langle a_n,b_k\rangle$ and only
  those $s$ with $h_s\in O$ respectively in $E$, are selected.

  Now since $g_0=g_1=1$, there must exist $s<t$ such that $c_s,c_t$
  are both $b_*^*$ or both $a_*^*$ and $h_s=h_t$, both of them
  occurring in the same $g_i$, and none of the $c_{s+1},\dots,c_{t-1}$
  contributing to that $g_i$. The relations in $\mathcal R$ and
  $\overline F$ allow the commutation of $c_s$ with
  $c_{s+1},\dots,c_{t-1}$, and its eventual cancellation with $c_t$.
  One then proceeds by induction on the length of $g$.
\end{proof}

\begin{theorem}
  The group $\kwv$ has the following presentation:
  \[
  \kwv=\Big\langle a_1,\ldots a_n,b_1,\dots,b_k\,\Big|\,\{\text{relations of }\overline F\}\cup\varphi^\ell(\rel)\text{ for all
  }\ell\ge0\Big\rangle,
  \]
  where $\varphi$ and $\rel$ and given respectively before
  Proposition~\eqref{prop:kwv:varphi} and in~\eqref{eq:kv:R}.
\end{theorem}
\begin{proof}
  This follows from Lemmata~\ref{lem:kwv:2}, \ref{lem:kwv:3}
  and~\ref{lem:kwv:4}.
\end{proof}

Since the endomorphism $\varphi$ is injective, we can embed the group
$\kwv$ into its ascending HNN-extension by $\varphi$, i.e.\ in the
group generated by $\kwv$ and an element $t$ whose action by
conjugation on $\kwv$ coincides with $\varphi$.

There is an $u\in\langle b_2,b_2^{b_1},b_3,b_3^{b_1}\rangle$ such that
$\varphi(b_k)=a_1^{ub_1^{y_k}}$; the group $\langle a_1^u,b_1\rangle$
is dihedral of order $2^{m+2}$, and we have $(a_1^ub_1)^{2^{m+1}}=1$.
We set $a=a_1^u$ and $b=b_1$ to obtain a new generating system
$\{b,b^t,\dots,b^{t^{k-1}}, a, a^t, \ldots, a^{t^{n-1}}\}$ of the
group $\kwv$; the identification with $\{a_i,b_j\}$ is similar to that
in~\eqref{eq:kv:identify}.  The last theorem then yields
\begin{theorem}
  Write $p(t)=x_n+x_{n-1}t+\dots+x_1t^{n-1}\in\Z[t]$ and
  $q(t)=y_k+y_{k-1}t+\dots+y_1t^{k-1}\in\Z[t]$. Let $m\in\{1,2\}$ and
  $u\in\langle b^t,b^{tb},b^{t^2},b^{t^2b}\rangle$ be as above. Then
  the group $\kwv$ is isomorphic to the subgroup $\langle a, a^t,
  \ldots, a^{t^{n-1}},b,b^t,\dots,b^{t^{k-1}}\rangle$ of the finitely
  presented group (where the generator $a$ is redundant)
  \begin{align*}
    \Big\langle a, b,t\,\Big|\,a^2,b^2,(ab)^{2^{m+1}},&
    b^{t^k}a^{b^{q(t)}},a^{u^{-1}(t^n-b^{p(t)})},
    [a,a^{t^{jn/d}}]\text{ for all }j\in\{1,\dots,d-1\},\\
    &\big[b^{t^i},b^{t^jb(ab)^{2\ell}}\big]
    \text{ for all }1\le i,j<k\text{ and }0\le\ell\le2^m,\\
    &\big[a^{t^i},b^{t^jb(ab)^{2\ell}}\big]
    \text{ for all }1\le i<n,1\le j<k\text{ and }0\le\ell\le2^m,\\
    &\big[a^{t^i},a^{t^jb(ab)^{2\ell}}\big]
    \text{ for all }1\le i,j<n\text{ and }0\le\ell\le2^m\Big\rangle.
  \end{align*}
\end{theorem}
\begin{proof}
  The first three relations express $\langle a_1^u,b_1\rangle$ as a
  dihedral group. The fourth one is $b_k^t=a_1^{ub_1^{y_k}}$, and the
  fifth one is $a_n^t=a_1^{b_1^{x_n}}$. The next ones are relations in
  $F$, and the last three rows are the commutation relations of the
  form $[a_i,a_j^w]$ and $[a_i,b_j^w]$ for $w\in E$ or $\in O$.
\end{proof}

\subsection{Branchness} We see in this subsection that $\kwv$ is
branch as soon as $k>1$ or $n>1$:
\begin{theorem}\begin{enumerate}
  \item If $k\ge 2$ and $n\ge2$, or if $k\ge3$ and $n=1$, then $\kwv$
    is branch on $\kwv'$; we have $\kwv/\kwv'\cong(\ZZ2)^{k+n}$
    generated by $\{a_i,b_j\}$ and
    $\kwv'/(\kwv'\times\kwv')=(\ZZ2)^{k+n-1}$ generated by
    $\{[b_1,b_j]_{2\le j\le k},[b_1,a_i]_{1\le i\le n}\}$.
  \item If $k=2$ and $n=1$, then $\kwv$ is branch on
    $L=\langle[b_1,b_2a_1]\rangle^{\kwv}$; we have $\kwv'/L=(\ZZ2)^2$
    generated by $\{[b_1,b_2],[b_2,a_1]\}$. Set $x=[b_1,b_2a_1]$; then
    $L/(L\times L)=\ZZ4\times\ZZ2$, generated by $\{x,x^{a_1}\}$ with the
    relations $x^4=(x^{a_1})^4=x^2x^{2a_1}=1$.
  \item If $k=1$ and $n\ge2$, then $\kwv$ is branch on
    $L=\langle[a_i,a_j]_{1\le i<j\le n},[a_i,b_1]_{1\le
      i<n}\rangle^{\kwv}$; we have $\kwv'/L=2$, generated by
    $[a_n,b_1]$, and $L/(L\times L)=\ZZ4\times(\ZZ2)^{n-2}$, generated by
    $\{[a_i,b_1]\}$ where $[a_1,b_1]$ has order $4$ and the other
    generators have order $2$.
  \item If $k=n=1$, then $\kwv$ is infinite dihedral, and is not even
    weakly branch.
  \end{enumerate}
\end{theorem}
\begin{proof}
  We use the endomorphism from Proposition~\ref{prop:kwv:varphi}. In all
  cases, it is easy to check that $L$ has finite index as claimed, and
  that $\varphi(L)\subset\kwv\times1$.

  It then remains to check that $\varphi(L)\subset L$. This is obvious
  if $L=\kwv'$ is characteristic. If $k=2$ and $n=1$, we obviously
  have $\varphi(\gamma_3(\kwv))\subset\gamma_3(\kwv)$; but we also
  have $\varphi(x)=[b_2^r,a_1^ta_1^u]\in\gamma_3(\kwv)$, so
  $\varphi(L)\subset L$ with $L=\langle\gamma_3(\kwv),x\rangle$ as
  claimed.

  If $k=1$ and $n\ge2$, then $\varphi[a_i,a_j]=[a_{i+1}^s,a_{j+1}^t]$
  for some $s,t\in\langle b_1\rangle$; now $i+1<n$ or $j+1<n$, so
  $[a_{i+1},a_{j+1}]\in L$ and either $[a_{i+1},b_1]\in L$ and
  $[a_{j+1},b_1]\in L$, so $[a_{i+1}^s,a_{j+1}^t]\in L$. Similarly
  $\varphi[a_i,b_1]=[a_{i+1}^s,a_1^t]\in L$.
\end{proof}

\subsection{Torsion}
First, we note that the group $\kwv$ always contains elements of
infinite order:
\begin{proposition}
  Every element $x$ which is a product, in any order, of all the
  generators $a_i,b_j$ of $\kwv$ has infinite order.
\end{proposition}
\begin{proof}
  Consider such an $x$. Its image in $\autxs/\autxs'$ is $(1,1,\dots)$
  so this element acts level-transitively; in particular, it has
  infinite order.
\end{proof}

\begin{proposition}
  If $k\ne1$ or $n\ne1$, then $\kwv$ contains every finite $2$-group
  as a subgroup. In particular, it contains torsion elements of
  arbitrarily large order.
\end{proposition}
\begin{proof}
  This follows directly from Proposition~\ref{prop:containsW}, since
  the groups $\kwv$ are regular branch on a subgroup $L$ containing a
  torsion element.

  If $k\ge2$ and $n\ge2$, this torsion element may be chosen as
  $[b_1,b_2]$. If $k=2$ and $n=1$, we may check that $(b_1a_1)^4$
  belongs to $L$ and has order $2$. If $k=1$ and $n\ge2$ then
  $[b_1,a_1]$ belongs to $L$ and has order $4$.
\end{proof}

\section{Kneading sequence and quadratic polynomials}\label{ss:kneading}
\subsection{Review of results in holomorphic dynamics}
Let $f(z)=z^2+c$ be a quadratic polynomial. Suppose that the orbit
of the critical point $0$ under the iterations of $f$ is finite.
Such polynomials are called \emph{post-critically finite}. We
distinguish two cases: when $0$ belongs to a finite cycle
(\emph{periodic case}) and when it does not, but its orbit is
still finite (\emph{pre-periodic case}).

Recall that if $M\subset\C$ is a connected and closed set with
connected complement, then there exists a unique biholomorphic
isomorphism $\Phi_M$ of the complement $\overline\C\setminus M$ with
$\{z\in\overline\C\;:\;|z|>R\}$ such that $\Phi_M(\infty)=\infty$ and
$\Phi_M'(\infty)=1$. The \emph{external ray} $\mathbf{R}_\alpha$ is the image
of the ray $\{r\cdot e^{2\pi i\cdot\alpha}\;:\;r\in (R, \infty)\}$
under $\Phi_M^{-1}$. One says that an external ray $\mathbf{R}_\alpha$
\emph{lands}, if the limit $\lim_{r\searrow R}\Phi_M^{-1}(r\cdot
e^{2\pi i\alpha})$ exists. Here and below the angle $\alpha$ is
considered to be an element of the group $\T$, i.e., the angles are
counted in full turns.

We use known facts about the dynamics of iterations of quadratic
maps (see~\cite{DH:orsayI,DH:orsayII,henkdierk}). The filled-in
Julia set $K_c$ of $z^2+c$ is the set of points which do not
escape to infinity under iteration, and the Fatou set is the open
set $\overline\C\setminus\partial K_c$. We assume that $z^2+c$ is
post-critically finite, hence $K_c$ is connected. External rays to
the Mandelbrot set are called \emph{parameter rays}, and external
rays to $K_c$ are called \emph{dynamical rays}.

It is easy to see that the image of a dynamical ray
$\mathbf{R}_\alpha$ under the action of $z^2+c$ is equal to the
ray $\mathbf{R}_{2\alpha}$.

Suppose that $0$ belongs to a cycle of length $n$ under iteration
of $z^2+c$. Then $c$ belongs to a \emph{hyperbolic component}
$M_c$ of the interior of the Mandelbrot set. For any other point
$c_1$ of that component, the quadratic polynomial $z^2+c_1$ also
has a unique attracting cycle of length $n$. If $\Phi(c_1)$
denotes the multiplier of this cycle (i.e., the product of
derivatives in all points of the cycle), then $\Phi$ is a
conformal isomorphism of $M_c$ with the open unit disc
$\mathbb{D}=\{z\in\C\;:\;|z|<1\}$. We obviously have
$\Phi(c)=0$, hence $c$ is called the \emph{center} of the
hyperbolic component $M_c$. The isomorphism $\Phi:M_c\to
\mathbb{D}$ extends to a homeomorphism of the boundary of $M_c$
with the unit circle. The preimage of $1$ under this homeomorphism
is called the \emph{root} of the component $M_c$. There exist
exactly two angles $\theta$ such that the parameter ray
$\mathbf{R}_\theta$ lands on the root of $M_c$.

In the dynamical plane, the point $c$ belongs to a Fatou component
$U_c$, which is periodic with period $n$ under $f$. There is a
unique point $r$ on the boundary of $U_c$, fixed under the map
$f^n:U_c\to U_c$ (since $f^n|_{U_c}$ is topologically conjugate
(via the \emph{B\"ottcher map}) to the restriction of $z^2$ to
$\mathbb{D}$). This point and its forward images are called the
\emph{roots} of their correspondent Fatou components.

A parameter ray $\mathbf{R}_\theta$ lands on the root of the
hyperbolic component $M_c$ if and only if the dynamical ray
$\mathbf{R}_\theta$ lands on the root of the Fatou component
$U_c$. Moreover, the number $\theta\in\T$ belongs to a cycle of
length $n$ under the doubling map $\alpha\mapsto 2\alpha:\T\to\T$.
In particular the angle $\theta$ is equal to $p/(2^n-1)$ for some
integers $p,n$, and the ray $\mathbf{R}_{2^k\theta}$ lands at the root of
the Fatou component to which $f^k(c)$ belongs.

Conversely, for every rational number $\theta\in\T$ with odd
denominator, the parameter ray $\mathbf{R}_\theta$ lands on the
root of a hyperbolic component $M_c$, and if $c$ is the center of
the component (i.e., the preimage of 0 under the multiplicator
map), then 0 has the same period under $z^2+c$ as has $\theta$
under the doubling map, and the dynamical ray $\mathbf{R}_\theta$
lands on the root of the Fatou component of $z^2+c$ containing
$c$.

Suppose now that $0$ is pre-periodic. Then $c$ belongs to the
boundary of the Mandelbrot set (it is a \emph{Misiurewicz} point) and
there exists a finite set of angles $\theta$ such that the
parameter rays $\mathbf{R}_\theta$ land on $c$. For each such
$\theta$ the external ray $\mathbf{R}_\theta$ in the dynamical
plane of $z^2+c$ lands on $c$. The pre-period of $\theta$ under
the doubling map is the same as the pre-period of $c$ under
$z^2+c$, but the period of $\theta$ may be a multiple of the
period of $c$. Here \emph{pre-period} and \emph{period} of a point
$x$ under a map $f$ are the minimal positive integers $k$ and $n$
such that $f^{k+n}(x)=f^k(x)$.

For example, the point $c\approx -0.1011+0.9563i$ is the landing
point of the parameter rays $\mathbf{R}_\alpha$ for
$\alpha=\frac{9}{56}, \frac{11}{56}$ and $\frac{15}{56}$. The
point $c$ has pre-period of length 3 and period of length 1 (i.e.,
it lands on a fixed point). But the angles have period 3, namely
\[9/56\mapsto 9/28\mapsto 9/14\mapsto 2/7\mapsto 4/7\mapsto 1/7\mapsto 2/7.\]

The period of $c$ is determined by $\theta$ as the period of the
\emph{kneading sequence} of $\theta$. Let $S_0$ be the image in
$\T$ of the interval $[\theta/2,(1+\theta)/2]$, and let $S_1$ be
the image in $\T$ of $[(1+\theta)/2,(2+\theta)/2]$. For every
$\alpha\in\T$ denote by $I_\theta(\alpha)$ its
\emph{$\theta$-itinerary}, defined as the sequence $a_0a_1\ldots$,
where
\[
a_k=\begin{cases} 0 & \text{ if }2^k\alpha\in S_0,\\
  1 & \text{ if }2^k\alpha\in S_1,\\
  * & \text{ if }2^k\alpha\in\{\theta/2, (1+\theta)/2\}.
\end{cases}
\]

The itinerary $I_\theta(\theta)$ is called the \emph{kneading
  sequence} of the point $\theta\in\T$ and is denoted
$\widehat\theta$.

If $\theta$ is periodic under the doubling map with period of
length $n$, then its kneading sequence is of the form $vv\ldots$,
where
\[
v=1x_2\ldots x_{n-1}{*},
\]
for some $x_i\in\{0, 1\}$.

If $\theta$ is strictly pre-periodic with a pre-period of length $k$
and period of length $n$, then its kneading sequence is of the form
$wvv\ldots$, where $w=1x_2\ldots x_k$ and $v$ are some words over
$\{0, 1\}$ of length $k$ and $n$ respectively with different last
letters.

It may happen that $v$ is a proper power. Then the period of the
kneading sequence is a factor of the period of $\theta$ under the
doubling map. In any case if $\theta$ is pre-periodic (i.e., if its
smallest denominator is even), and the parameter ray
$\mathbf{R}_\theta$ lands on $c$, then the period of $c$ under
iteration of $z^2+c$ is equal to the period of the kneading sequence
$\widehat\theta$.

\subsection{Iterated monodromy groups of quadratic polynomials}
Suppose that $w$ is either a kneading sequence of the form
$(x_1\ldots x_{n-1}*)^\omega$, or a kneading sequence of the form
$y_1\ldots y_k(x_1\ldots x_n)^\omega$, where $y_k\ne x_n$ and
$x_1\ldots x_n$ is not periodic, i.e., is not a proper power (note
that every pre-periodic sequence can be uniquely represented in
that form). Then we denote by $\kn{w}$ the group $\kn{x_1\ldots
x_{n-1}}$ in the first case and $\kn{y_1\ldots y_k,\,x_1\ldots
x_n}$ in the second.

\begin{theorem}\label{th:img}
  Let $f(z)=z^2+c$ be a post-critically finite quadratic polynomial.
  Let $\theta\in\T$ be an angle such that the parameter ray
  $\mathbf{R}_\theta$ lands either on the root of the hyperbolic
  component $M_c$ (if $c$ is periodic) or on $c$ (if $c$ is
  pre-periodic).

  Then $\img{z^2+c}$ is isomorphic to $\kn{\widehat{\theta}}$.
  Moreover, the action of $\img{z^2+c}$ on the tree of preimages is
  conjugate with the action of $\kn{\widehat{\theta}}$ on the binary
  tree.
\end{theorem}
\noindent We consider independently the periodic and pre-periodic cases.
\begin{proof}[Proof in the periodic case]
  We use the \emph{invariant spiders}, described
  in~\cite{hubbardschleicher}, to make cuts in $\overline\C$. A
  \emph{spider} is a collection of disjoint closed paths $\gamma_z$,
  called \emph{legs}, connecting every point $z\in P_f$ to infinity. A
  spider $\mathcal{S}$ is $f$-invariant if
  $f^{-1}(\mathcal{S})\supset\mathcal{S}$, up to an isotopy
  relative to $P_f$.

  The dynamical ray $\mathbf{R}_\theta$ has two preimages under $f$:
  the ray $\mathbf{R}_{\theta/2}$ and $\mathbf{R}_{(1+\theta)/2}$.
  Both rays land in the pre-periodic case on $0$ and divide the plane
  into two connected components.

  In the periodic case, the rays $\mathbf{R}_{\theta/2}$ and
  $\mathbf{R}_{(1+\theta)/2}$ land on two points belonging to the
  boundary of the Fatou component containing $0$. These two points are
  the preimages of the root of the component of $c$. Let us connect
  zero to these two points by \emph{internal rays} of the Fatou
  component, i.e.\ by images of rays under the holomorphic map fixing
  $0$ and mapping the component to the disc $\mathbb{D}$.

  The union of the rays $\mathbf{R}_{\theta/2}$,
  $\mathbf{R}_{(1+\theta)/2}$ and the constructed internal rays also
  divide the plane into two connected components. Let us denote the
  component to which $c$ belongs by $S_1$, and the other component (to
  which the ray $\mathbf{R}_0$ belongs) by $S_0$.

  For each $i\in\{0,1\}$, the restriction of $z^2+c$ to the component
  $S_i$ is a homeomorphism of $S_i$ onto the set
  $\C\setminus\gamma_c$, where $\gamma_c$ is a curve
  connecting $c$ to infinity. The curve $\gamma_c$ is the union of the
  dynamical ray $\mathbf{R}_\theta$ with an internal ray.

  Let $t\notin\bigcup_{n\ge 0}f^k(\gamma_c)$ be a basepoint and let
  $T=\bigsqcup_{n\ge 0}f^{-n}(t)$ be the corresponding tree of
  preimages. We will construct an isomorphism between the tree $T$ and
  the binary tree $\{0, 1\}^*$ using \emph{itineraries}.

  If $z\in f^{-n}(t)$ is a vertex of the tree $T$, then the
  corresponding vertex of the binary tree is given by the word
  $\Lambda(z)=x_{n-1}x_{n-2}\ldots x_0$, where $x_k\in\{0, 1\}$ is
  such that
  \[f^k(z)\in S_{x_k}.\]

  It follows directly from the definition that $\Lambda$ is a
  level-preserving bijection and that if $\Lambda(z)=x_{n-1}\ldots
  x_0$, then $\Lambda(f(z))=x_{n-1}\ldots x_1$. Hence
  $\Lambda:T\arr\{0, 1\}^*$ is an isomorphism of rooted trees.  From
  now on, we identify the trees $\xs$ and $T$ by this isomorphism.

  Let $P_f=\bigcup_{i\ge 1}f^r(c)$ be the post-critical orbit of $f$,
  and set $\M=\C\setminus P_f$, so that $f$ is a covering map
  from $f^{-1}(\M)$ to $\M$. If $n$ is the period of the point $c$,
  then $n$ is also the period of the curve $\gamma_c$ under iteration
  of $f$. Moreover, the curves $f^k(\gamma_c)$ are pairwise disjoint
  for $k=0, 1, \ldots, n-1$. More precisely, they can have common
  points, but they do not intersect transversally, so that they become
  pairwise disjoint after small homotopies in $\M$,
  see~\cite{hubbardschleicher}. Then the set $\{f^k(\gamma_c)\;:\;0\le
  k\le n-1\}$ is an invariant spider. For every $z\in P_f$ we denote
  by $\gamma_z$ the unique leg $f^k(\gamma_c)$ of the spider which
  connects the point $z$ to infinity. We have then
  $f(\gamma_z)=\gamma_{f(z)}$. The right-hand side of
  Figure~\ref{fig:spider} shows the preimage of a spider.

  An example of spider $\bigcup_{z\in P_f}\gamma_z$ is shown on the
  left hand side of Figure~\ref{fig:spider}, for
  $c\approx-0.1225+0.7448i$ (the ``Douaddy rabbit'').  The
  corresponding $\theta$ is equal either to $1/7$ or to $2/7$ (it is
  $2/7$ in our picture), the critical point belongs to a cycle of
  length $3$ and the root of all three components is a common fixed
  point.

  \begin{figure}[ht]
    \begin{center}
      \includegraphics{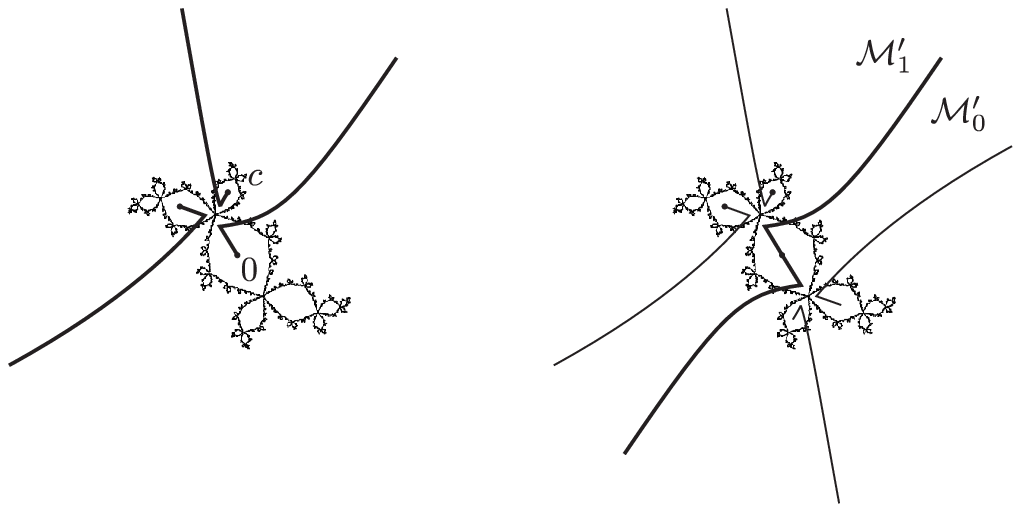}
    \end{center}
    \caption{Invariant spider}\label{fig:spider}
  \end{figure}

  Write $\M'=\M\setminus\bigcup_{z\in P_f}\gamma_z$. It follows that
  $\M'$ is simply connected and $f^{-1}(\M')\subset\M'$, up to
  isotopy. The set $f^{-1}(\M')$ has two connected components
  $\M'_0=S_0\cap\M'$ and $\M_1'=S_1\cap\M'$.

  Define for $t, u\in\M'$ and $z\in P_f$ the path $g_z(t,u)$ starting
  in $t$, ending in $u$, disjoint from all the paths $\gamma_y$ of the
  spider for $y\ne z$, and intersecting $\gamma_z$ only once in such a
  direction that the part of $\gamma_z$ containing $z$ is to the left
  of the path $g_z(t,u)$. The path $g_z(t,u)$ is obviously defined
  uniquely up to homotopy in $\M$. We also denote by $g(t,u)$ the path
  starting in $t$, ending in $u$ and disjoint from the legs of the
  spider. The path $g(t,u)$ is also uniquely defined up to homotopy in
  $\M$. If $y\notin P_f$ then $g_y(t,u)$ is, by definition, the path
  $g(t,u)$.

  It follows from the definitions that, for $z\ne c$,
  \[f^{-1}(g_z(t,u))=\{g_{z_0}(t_0,u_0), g_{z_1}(t_1,u_1)\},\] where
  $f^{-1}(z)=\{z_0,z_1\}$ and $f^{-1}(t)=\{t_0,t_1\}$ and
  $f^{-1}(u)=\{u_0,u_1\}$. We may assume, up to changing our notation,
  that $z_0\in S_0$ and $z_1\in S_1$. Then $t_0,u_0\in S_0$ and
  $t_1,u_1\in S_1$.

  If $z=c$, then the only preimage of $z$ is $0$ and only one of the
  components of $f^{-1}(g_c(t,u))$ intersects the path $\gamma_0$.
  Then this component is homotopic to $g_0(t_0,u_0)$, and the other
  component is homotopic to $g(t_1,u_1)$.  If $\gamma_0$ is the
  extension of the ray $\mathbf{R}_{\theta/2}$, then $t_0,u_1\in S_0$
  and $t_1,u_0\in S_1$. If $\gamma_0$ is the extension of the ray
  $\mathbf{R}_{(1+\theta)/2}$, then $t_0,u_1\in S_1$ and $t_1,u_0\in
  S_0$.

  We also have that the $f$-preimages of the path $g(t,u)$ are the
  paths $g(t_0, u_0)$ and $g(t_1,u_1)$, where $t_0,u_0\in S_0$ and
  $t_1,u_1\in S_1$.

  It is not hard to see now that if $\gamma_0$ is the extension of
  the ray $\mathbf{R}_{(1+\theta)/2}$, then the generators
  $h_k=g_{f^k(0)}(t,t)$ of $\pi_1(\M, t)$ act on the tree $\xs$ in the
  same way as the generators $a_k$ of the group
  $\kn{\widehat{\theta}}$. If $\gamma_0$ is the extension of the ray
  $\mathbf{R}_{\theta/2}$, then the generators $h_k$ act in the same
  way as the generators $a_k^{-1}$.

  We conclude that the actions of $\img{f}$ on $T$ and of
  $\kn{\widehat{f}}$ on $\xs$ are conjugate. Note that the proof
  does not depend on the choice of $\theta$ for a given $c$
  (though the curve $\gamma_0$ does depend on $\theta$). The
  wreath recursion defining $\kn{\widehat{\theta}}$ is given by
  the paths disjoint with the legs of the spider connecting the
  basepoint to its preimages. A different choice of the
  connecting paths (and hence of the wreath recursion) are
  convenient in some other situations (see, for
  instance~\cite{bartnek:rabbit}).
\end{proof}

\begin{proof}[Proof in the pre-periodic case]
  The problem here is that there is no invariant spider with disjoint
  legs, when the period of the angle is greater than that of the
  kneading sequence.

  However, we can find a sequence of spiders $\mathcal{S}_0,
  \mathcal{S}_1, \ldots$, such that $\mathcal{S}_{k+1}\subset
  f^{-1}(\mathcal{S}_k)$. Take any spider $\mathcal{S}_0$
  consisting of dynamical rays landing on $P_f$ and define inductively
  $\mathcal{S}_k$ to be the set of paths belonging to
  $f^{-1}(\mathcal{S}_{k-1})$ and landing on the points of
  $P_f$. The points of $P_f$ are not critical, therefore for every
  $z\in P_f$ there exists a unique path $\gamma_{z, k+1}\in
  f^{-1}(\mathcal{S}_{k-1})$ landing on $z$. Hence the
  spider $\mathcal{S}_k$ is well defined.

  We can then define the paths $g_{z, k}(t,u)$ for every spider
  $\mathcal{S}_k$ in the same way as in the periodic case. If
  $\gamma_{c, k}$ is the element of $\mathcal{S}_k$ landing on $c$,
  then we denote by $S_{0, k+1}$ and $S_{1, k+1}$ the components of
  $\C\setminus f^{-1}(\gamma_{c, k})$, where $S_{1,
    k+1}$ is the component containing $c$. The isomorphism
  $\Lambda:T\arr\xs$ is then defined using the itineraries of points
  with respect to these partitions of the plane. Namely, if $z\in
  f^{-k}(t)$ is a vertex of the tree $T$, then
  $\Lambda(z)=x_{n-1}x_{n-2}\ldots x_0$, where $x_k\in\{0, 1\}$ is
  such that
  \[f^k(z)\in S_{x_k, n-k}.\] It is also easy to prove that $\Lambda$
  is an isomorphism of rooted trees.

  The same formulae for the preimages of the paths $g_{z, k}(t,u)$
  hold as in the periodic case. The only differences will be that the
  index $k+1$ appears at the preimages of paths and at names of the
  components $S_{i, k}$, and that the point $0$ does not belong to the
  post-critical set. Therefore the preimages of $g_c(t,u)$ are the
  paths $g(t_0,u_0)$ and $g(t_1,u_1)$, with $t_0,u_1\in S_{0, k+1}$
  and $t_1,u_0\in S_{1, k+1}$.  The partitions of the plane into
  components $S_{0, k+1}$ and $S_{1, k+1}$ agree with the kneading
  sequence (i.e., every $z\in P_f$ belongs either only to the sectors
  $S_{0, k}$ or only to the sectors $S_{1, k}$), since the kneading
  sequences of all rays landing on $c$ are equal to
  $\widehat{\theta}$. These considerations prove that the generators
  $g_{z, 0}(t, t)$ of $\img{f}$ act in the same way as the generators
  $a_i, b_i$ of the group $\kn{\widehat{\theta}}$.
\end{proof}

The following result follows from the classical results in symbolic
dynamics of quadratic polynomials (see~\cite{henkdierk,keller}). It
follows also from general results on iterated monodromy groups of
expanding maps (see~\cite[Theorem~6.4.4]{nek:book}).

\begin{theorem}
  Under the conditions of Theorem~\ref{th:img}, the limit dynamical
  system of the group $\kn{\widehat{\theta}}$ is topologically
  conjugate to the action of the polynomial $z^2+c$ on its Julia set.
\end{theorem}

\bibliographystyle{plain}
\bibliography{mymath,nekrash}
\end{document}